\documentclass{amsproc}
\usepackage{amssymb, epsf}
\newtheorem{theorem}{Theorem}             
\newtheorem{proposition}[theorem]{Proposition}
\newtheorem{lemma}[theorem]{Lemma}
\newtheorem{corollary}[theorem]{Corollary}
\theoremstyle{remark}

\newtheorem{example}{Example}

\newcommand\lbb[1]{\label{#1}}% {\texttt{(#1)}}\;\;\;}    % temporary
\newcommand\ov[1]{\overline{#1}}

\newcommand{\fig}[1]
        {\raisebox{-0.5\height}%
                  {\epsfbox{#1.eps}}%
        }
\def\H{\mathcal H}      % Hecke algebra
\def\Z{{\mathbf Z}}     % integers
\def\l{\lambda}
\def\eps{\varepsilon}
\def\si{\sigma}
\def\vac{{\mathbf 1}}
\def\vi{v^{-1}}
\def\Zv{{\mathbf Z[v^{\pm 1}]}}  % 
\def\vzv{{ v\mathbf Z[v]}}  % 
\def\Qv{{\mathbf Q(v)}}          % 
\def\KL{ Kazhdan-Lusztig }
\def\P{\mathcal P}  
\def\d{\partial}
\def\s{\scriptstyle }
\begin{document}
\title{ Factorization of Kazhdan--Lusztig elements for Grassmanians}
\author[A. Kirillov, Jr. and A. Lascoux]{Alexander Kirillov, Jr. \and
          Alain Lascoux} 
\address{
{\rm Alexander Kirillov} \\
Institute for Advanced Study, Princeton, NJ 08540, USA}
\email{kirillov@math.ias.edu}
\address{
{\rm Alain Lascoux}\\
CNRS, Institut Gaspard Monge, Universit\'e de Marne-la-Vall\'ee, \\
         5, boulevard Descartes, 77454 Marne-la-Vall\'ee, Cedex 2, France}
\email{Alain.Lascoux@univ-mlv.fr}
\thanks{The first author was partially supported by NSF grants
  DMS-9610201, DMS-97-29992}

\begin{abstract}
We show that the Kazhdan-Lusztig basis elements $C_w$ of the Hecke
algebra of the symmetric group, when $w \in S_n$ corresponds to a
Schubert subvariety of a Grassmann variety, can be
written as a product of factors of the form $T_i+f_j(v)$, where $f_j$
are rational functions.
\end{abstract}
\maketitle

\section{Notation} \lbb{snota}
In this section, we briefly list the main facts and notations related
to  Kazhdan--Lusztig polynomials
and their parabolic analogues  (see \cite{D}, \cite{S}). We use the following
notations:

$\H$---the Hecke algebra of the symmetric group $S_n$; we consider it as
an algebra over the field $\Qv$ (the variable $v$ is related to the
variable $q$ used by Kazhdan and Lusztig via $v=q^{1/2}$),
and we write the quadratic relation in the form 
\begin{equation*}
(T_i-v)(T_i+\vi)=0.
\end{equation*}

$C_w$---KL basis in $\H$, which we define by the conditions
$\ov{C_w}=C_w$, $C_w-T_w\in \oplus \vzv T_y$.

For any subset $J\subset \{1, \dots, n-1\}$, we denote by $W_J\subset
S_n$ the corresponding parabolic subgroup, and by $W^J$ the set of
minimal length representatives of cosets $S_n/W_J$. We also denote by
$M^J$ the $\H$-module induced from the one-dimensional representation
of $\H(W_J)$, given by $T_j m_1 =-\vi m_1, j\in J$. We denote
$m_y=T_y m_1, y\in W^J$ the usual basis in $M^J$. 

We define the parabolic KL basis $C_y^J, y\in W^J$ in $M^J$ by
$\ov{C^J_y}=C^J_y, C^J_y-m_y\in \oplus_{z\in W^J} \vzv m_z$. 

Denote for brevity $C_J=C_{w_0^J}$ the element of KL basis in $\H$
corresponding to the element of $w_0^J$ of maximal length in
$W_J$. The following result is well-known (see, e.g., \cite{S}).  

\begin{lemma}\lbb{l1}

\textup{(i)}
\begin{equation*}
C_J=\sum_{w\in W^J}(-v)^{l(w_0^J)-l(w)}T_w.
\end{equation*}

\textup{(ii)} Let $w\in W$ be such that it is an element of maximal
length in the coset $w W_J$ \textup{(}which is equivalent to $w=\tau
w_0^J$ for some $\tau\in W^J$\textup{)}. Then $C_w=X C_J$ for some
$X\in \oplus_{y\in W^J} \Zv T_y$.

\textup{(iii)} Let $X\in \oplus_{y\in W^J} \Zv T_y$. Then 
\begin{equation*}
X m_1=C^J_\tau \iff XC_J=C_{\tau w_0^J}.
\end{equation*}
\end{lemma}

Let us now consider the special case of the above situation. From now
on, fix $k\le n-1$, and let $J=\{1, \dots, k-1, k+1, \dots, n-1\}$ so
that $W_J=S_k\times S_{n-k}$ is a  maximal parabolic subgroup in
$S_n$. In this case, the module $M^J$ can be described as follows: 
\begin{equation}
\begin{aligned}
M=&\bigoplus_{\eps\in E} \Qv \eps,\\
T_i\eps&=\begin{cases} s_i \eps, 
                  & (\eps_i, \eps_{i+1})=(+-),\\
               -\vi \eps,
                  & (\eps_i, \eps_{i+1})=(--)\text{ or } (++),\\
               s_i\eps + (v-\vi) \eps,
                  &  (\eps_i, \eps_{i+1})=(-+),
          \end{cases}        
\end{aligned}
\end{equation}
where $E$ is the set of all length $n$ sequences of pluses and minuses
which contain exactly $k$ pluses. The relation of this with the
previous notation is given by $m_y\leftrightarrow
y(\vac )=T_y(\vac )$, where
\begin{equation}
\vac =(\underbrace{+\dots +}_{k}\underbrace{-\dots -}_{n-k}).
\end{equation}

In particular, $m_1\leftrightarrow \vac $. 

The set of minimal length representatives $W^J$ also admits a
description in terms of Young diagrams. Namely, let $\l$ be a Young
diagram which fits inside the $k\times (n-k)$ rectangle. Define
$w_\l\in S_n$
by 
\begin{equation}\lbb{wl}
w_\l=\prod_{(i,j)\in \l} s_{k+j-i},
\end{equation}
where $(i,j)$ stands for the box in the $i$-th row and $j$-th column,
and the product is taken in the following order: we start with the
lower right corner and continue along the row, until we get to the
first column; then we repeat the same with the next row, and so on
until we reach the upper left corner. 

\begin{example}\lbb{ex1}
Let $\l$ be the diagram shown below, and $k=7$ (to assist the reader,
we put the numbers $k+j-i$ in the diagram). 

\begin{equation*}
\fig{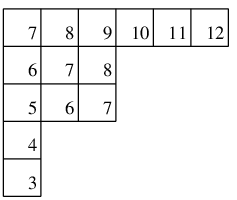}
\end{equation*}

Then 
$w_\l=s_3\cdot  s_4\cdot  s_7 s_6 s_5\cdot s_8s_7 s_6\cdot s_{12} s_{11} s_{10} s_9 s_8
s_7$ (for easier reading, we separated products corresponding to
different rows by $\cdot$). 
\end{example}

The proof of the following proposition is straightforward.

\begin{proposition}The corespondence $\l\mapsto w_\l$, where $w_\l$ is
  defined by \eqref{wl}, is a bijection between the set of all Young
  diagrams which fit inside the $k\times (n-k)$ rectangle  and $W^J$.  
\end{proposition}

\section{The main theorem}\lbb{smain}
As before, we fix $k\le n-1$ and let  $J=\{1, \dots, k-1, k+1, \dots,
n-1\}$. Unless otherwise specified, we only use Young diagrams which
fit inside the $k\times(n-k)$ rectangle.  

For a Young diagram $\l$, we define the shifts $r_{i,j}\in \Z_{>0},
(i,j)\in \l$ by the following relation
\begin{equation}\lbb{rij}
r_{ij}=\max (r_{i, j+1}, r_{i+1, j})+1,
\end{equation}
where we let $r_{ij}=0$ if $(i,j)\notin \l$. 
\begin{example}
For the diagram $\l$ from Example~\ref{ex1}, the shifts $r_{ij}$ are
shown below. 

\begin{equation*}
\fig{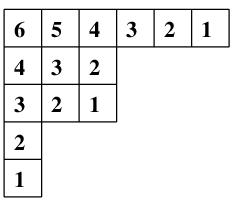}
\end{equation*}

\end{example}

Next, let us define for each diagram $\l$ an element $X_\l\in \H$ by 
\begin{equation}\lbb{xl}
X_\l=\prod_{(i,j)\in \l} \biggl( T_{k+j-i}
                                 -\frac{v^{r_{ij}}}{[r_{ij}]}
                         \biggr)
\end{equation}
where, as usual, $[r]=(v^r-v^{-r})/(v-\vi)$, and the product is taken
in the same order as in \eqref{wl}. 

The main result of this paper is the following theorem.
\begin{theorem} \lbb{main}
Let $\l$ be a Young diagram. Then 
\begin{equation*}
X_\l \vac =C^J_{w_\l}.
\end{equation*}
\end{theorem}
Note that by Lemma~\ref{l1}, this is equivalent to 
\begin{equation}
X_\l C_J=C_{w_\l w_0^J}.
\end{equation}

We remind the reader that the Kazhdan-Lusztig elements $C_{w w_0^J}$,
where $w\in W^J$, and $W_J$ is a maximal parabolic in $S_n$ (they are
also known as KL elements for Grassmanians), have been studied in a
number of papers. A combinatorial description was given in \cite{LS};
it was interpreted geometrically in \cite{Z}, and in terms of
representations of quantum $\mathfrak{gl}_m$ in \cite{FKK}. However,
it is unclear how these results are related with the factorization
given by the theorem above. A similar factorization was given in
\cite{ICM} for those permutations which correspond to non-singular
Schubert varieties---i.e., for those $w$ such that, for any $v\in
S_n$, the Kazhdan-Lusztig polynomial $P_{v,w}$ is either 1 or 0.

Note that one can easily check that the elements $X_\l$ are invariant
under the Kazhdan--Lusztig involution: $\ov{X_\l}=X_\l$; thus, all the
difficulty is in proving that they are integral and have the right
specialization at $v=0$.

A crucial step in proving this theorem is the following
proposition. 
\begin{proposition}\lbb{rect}
Theorem~{\rm\ref{main}} holds when $\l$ is the $k\times (n-k)$
rectangle.
\end{proposition}
\begin{proof}
For any $w\in S_n$, choose a reduced expression $w=s_{i_\ell}\dots
s_{i_1}$. Define the element $\nabla_w\in
\H$ by
\begin{equation}\lbb{nabla}
\nabla_w=\biggl(T_{i_\ell}-\frac{v^{r_\ell}}{[r_\ell]}\biggr)
         \dots 
         (T_{i_1}-v),
\end{equation}
where $r_1,\dots, r_\ell\in \Z_+$ are defined as follows: if
$s_{i_{m-1}}\dots s_{i_1}(1,\dots, n)=(\dots, a,b,\dots )$ (in
$i_{m}$-th, $(i_m+1)$-st places), then $r_m=b-a$. 
Then $\{\nabla_w, w\in S_n\}$ is a Yang-Baxter basis of the Hecke
algebra, and we have (see \cite[\S 3]{YB}): 
\begin{lemma}\lbb{yb}

\textup{(i)} The element $\nabla_w$ does not depend on the choice of
reduced expression. 

\textup{(ii)} If $w_0^J$ is the longest element in some parabolic
subgroup $W_J\subset S_n$, then $\nabla_{w_0^J}=C_J$. 
\end{lemma}

Now, let us prove our proposition, i.e. that $X_\l C_J$ is a KL
element for rectangular $\l$. In this case, $w_\l$ is the
longest element in $W^J$: 
\begin{equation*}
w_\l(\vac )=(\underbrace{-\dots -}_{n-k}\underbrace{+\dots +}_{k}).
\end{equation*}

Let us choose the following reduced expression for the longest element
$w_0$ in $S_n$: $w_0=w_\l w_0^J$, where we take for $w_\l$ the reduced
expression given by \eqref{wl}. Then one easily sees that definition
\eqref{nabla} in this case gives
\begin{equation*}
\nabla_{w_0}=X_\l\nabla_{w_0^J}.
\end{equation*}
By Lemma~\ref{yb}, we get $C_{w_0}=X_\l C_J$, which is exactly the
statement of the proposition. 
\end{proof}

The proof in the general case is based on the following
proposition. Denote 
\begin{equation}\lbb{ov}
O(v^m)=\{f\in \Qv | f \text{ has zero of order $\ge m$ at }v=0\}.
\end{equation}

\begin{proposition} \lbb{pov}
\begin{equation*}
X_\l \vac =w_\l(\vac )+\sum_{\eps\in E} O(v)\eps.
\end{equation*}
\end{proposition}
A proof of this proposition is given in Section~3. 

Now we can give a proof of the main theorem. First, one easily checks
the invariance under the bar involution, since
\begin{equation*}
\ov{T_i-\frac{v^r}{[r]}}=T_i-\frac{v^r}{[r]}.
\end{equation*}
Combining this with Proposition~\ref{pov}, we see that it remains to
show that $X_\l C_J$ are integral, i.e.  $X_\l C_J\in \oplus \Zv T_w$
(note that it is not true that $X_\l$ itself is integral.)
This will be done by induction.

Let $\l$ be a Young diagram. Then we claim that any  such diagram
can be presented as a union $\l=\l'\sqcup \mu$, where $\mu$ is a
rectangle, and $\l'$ is again a Young diagram such that for $(i,j)\in
\l'$, the shifts $r_{(i,j)}^{\l'}=r_{(i,j)}^{\l}$. It can be formally
proved as follows: if one 
writes the successive widths and heights of the stairs of the diagram 
$$ \infty ,(a_1,b_1),  (a_2,b_2), ...(a_k,b_k), \infty$$  
then there is at least one index $i$ for which 
$a_i \leq b_{i+1}$ and $b_i \leq a_{i+1}$.  In that case, the
rectangle $\mu$  has the lower right corner $i$.

\begin{example}
  For the diagram $\l$ from Example~\ref{ex1}, the sequence $(a_k,
 b_k)$ is given by $\infty, (1,2), (2,2), (3,1), \infty$, and the
 subdiagram $\mu$ is the shaded $2\times 2$ square, as shown below. As
 before, we also included the shifts $r_{ij}$ in this diagram. The
 subsets $I^\mu, J^\mu$ in this case are given by $I^\mu=\{6,7,8\},
 J^\mu=\{6,8\}$.
\begin{equation*}
\fig{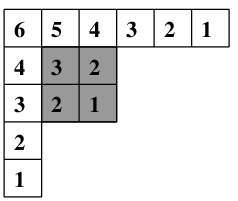}
\end{equation*}
\end{example}

Let us choose for $\l$ the presentation $\l=\l'\sqcup \mu$, where
$\mu$ is a rectangle, as above. Then $X_\l=X_\mu X_{\l'}$. 

Define the subsets $I^\mu, J^\mu\subset \{1, \dots, n-1\}$ by
$I^\mu=\{k'-a+1, \dots, k'+b-1\}, J^\mu=I^\mu\setminus\{k'\}$, where
$k'=k-i+j$, $(i,j)$---coordinates of the UL corner of $\mu$, 
$a$ and $b$ are numbers of rows and columns in $\mu$ respectively.

We need to show that $X_\mu X_{\l'}C_J\in \sum \Zv T_y$.  By induction
assumption, we may assume that $X_{\l'} C_J=C_{\si}$, where we denoted
for brevity $\si=w_{\l'}w_0^J$. It is easy to show that if $\mu$ is
chosen as before,  then $\si$ is the maximal length element in the
coset $W_{J^\mu}\si$. Thus, by Lemma~\ref{l1}, we can write $C_\si=
C_{J^\mu} Y$ for some integral $Y\in \H$. Therefore, 
$X_\mu X_{\l'}C_J=X_\mu C_{J^\mu}Y$. Since $W_{I^\mu}$ is itself a
symmetric group, and $W_{J_\mu}$ is a maximal parabolic subgroup in
it, we can use Proposition~\ref{rect}, which gives  $X_\mu
C_{J^\mu}=C_{I^\mu}$, and therefore, $X_\mu X_{\l'}C_J=C_{I^\mu}Y\in
\sum \Zv T_w$. 
\qed

\section{Proof of regularity at $v=0$}

In this section we give the proof of Proposition~\ref{pov}. Before
doing so, let us introduce some notation.

As before, assume that we are given $n, k, \l$ and a collection of
positive integers $r_{ij}, (i,j)\in \l$ (not necessarily defined as in
\eqref{rij}). Let $\eps\in E$ be a sequence of pluses and minuses. 
We define the {\it weight} $r_\l(\eps)$ as folllows.

Define $a(i), i=1\dots k$ by $a(i)=k+\l_i-i+1$. Equivalently, these
numbers can be characterized by saying that $w_\l(\vac )$ has pluses
exactly at positions $a(k), \dots, a(1)$. 

 Define $r_\l(\eps)=\sum_{t=1}^n r_t(\eps)$, where $r_t(\eps)$ is
   defined as 
   follows: 

(i) if $t=a(i), \eps_t=-$ then  $r_t(\eps)=r_{i,\l_i}-1$

(ii) if $a(i)<t<a(i+1), \eps_t=+$ then $r_t(\eps)=r_{i,j}, k+j-i=t$

(iii) otherwise, $r_t(\eps)=0$

In a sense, $r_\l(\eps)$ measures the discrepancy between $\eps$ and
$w_\l(\vac )$. Indeed, let us denote the numbers of rows and columns
in $\l$ by $i, j$ respectively, and let $\eps$ be such that
\begin{equation}\lbb{cond}
\begin{aligned}
\eps_t&=+ \text{ for } t\le k-i, \\
\eps_t&=- \text{ for } t> k+j. 
\end{aligned}
\end{equation}
 
Then one easily sees that
\begin{align}\lbb{r>0}
&r_\l(\eps)\ge 0,
&r_\l(\eps)=0\iff \eps=w_\l(\vac )
\end{align}

\begin{example}
Below we illustrate the calculation of $r_\l(\eps)$, where $\l$ is the
diagram used in Example~\ref{ex1}. The positions $a(i)$ are shaded
(thus, the sequence of colors encodes $w_\l(\vac )$, with
``shaded''$\leftrightarrow +$, ``unshaded''$\leftrightarrow -$), and
we connected unshaded pluses with the corresponding box $(i,j)$,
defined in (ii) above. For convenience of the reader, we also put the
numbers $k+j-i$ (not the shifts $r_{ij}$!) in the diagram.
\begin{equation*}
\fig{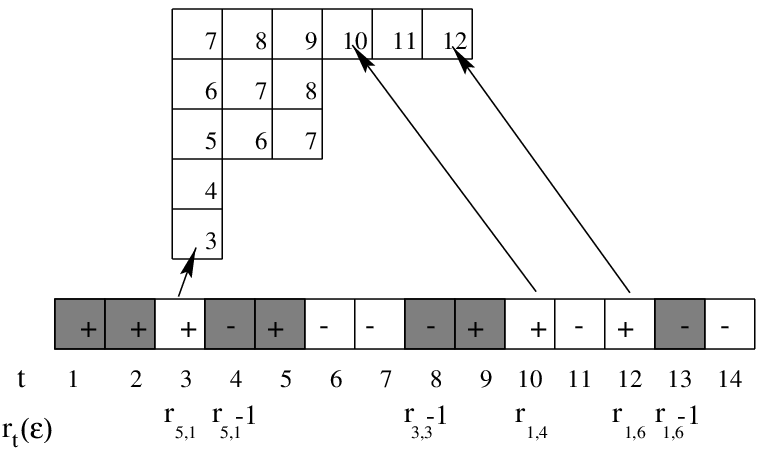}
\end{equation*}

\end{example}

\begin{lemma}Let $\l$ be any Young diagram inside the $k\times (n-k)$
  rectangle, and let $r_{ij},(i,j)\in\l$, be positive integers
  satisfying $r_{ij}>r_{i,j+1}, r_{ij}>r_{i+1, j}$. Define
  $\mathcal L_\l\subset M^J$ by  
\begin{equation*}
\mathcal L_\l=\sum_{\eps\in E}O(v^{r_\l(\eps)})\eps.
\end{equation*}
 Then
\begin{equation*}
X_\l \vac \in \mathcal L_\l.
\end{equation*}

\end{lemma}

Before proving this lemma note that due to \eqref{r>0},  
this lemma immediately implies
Proposition~\ref{pov}.

\begin{proof}
The proof is by induction. Let $(i,j)$ be a corner of $\l$, and
$\l'=\l - (i,j)$, so that 
$X_\l= \biggl(T_{k-i+j}-\frac{v^{r_{ij}}}{[r_{ij}]}\biggr)
        X_{\l'}$. 
Since $\frac{v^r}{[r]}\in O(v^{2r-1})$, it suffices to
prove  that 
$\biggl(T_{k-i+j}+O(v^{2r_{ij}-1})\biggr)\mathcal L_{\l'} 
\subset \mathcal L_{\l}$. Since this operation only changes $\eps_a,
\eps_{a+1}$ ($a=k-i+j$), we need to consider 4 cases: $(++), (+-),
(-+), (--)$. This is done explicitly. For example, for the $(+-)$
case, we have 
\begin{equation*}
\bigl(T_m+O(v^{2r_{ij}-1})\bigr)(\dots +-\dots )= 
(\dots -+\dots ) +O(v^{2r_{ij}-1})(\dots +-\dots )
\end{equation*}

In this case, the first  summand has the same weight and comes with
the same power of $v$ as the original $\eps$ (note that in the
original $\eps$, this $(+-)$ didn't contribute to the weight), so it
is in $\mathcal L_\l$. As for the second summand, its weight is increased by
$2r_{ij}-1$ (the plus contributes $r$ and the minus, $r-1$), but it
comes with the  factor $O(v^{2r_{ij}-1})$, so again, it is in
$\mathcal L_\l$. The other cases are treated similarly. 

\end{proof}

%%%%%%%%%%%%%%%%%%%%%%%%%%%%%%%%%%%%%%%%%%%%%%%%%%%%%%%%%%%%%%%%%%%%%%
\section{Divided differences and parabolic \KL bases}

In this section, we give a factorization for the dual Kazhdan--Lusztig
basis for Grassmanians. 

To induce a parabolic module, one can start from the 1-dimensional
representation $T_j\mapsto v$ instead of $T_j \mapsto -1/v$ which was
used in \S 1. 
We now denote the corresponding module by  $M'$ and its \KL basis by 
$C^{\prime J}_y$ to distinguish from previous case. Note that there exists a
natural pairing between $M$ and $M'$, and  $C^J_y$
and $C^{\prime J}_y$ are dual bases with respect to this pairing (see, e.g.,
\cite{S}, \cite{FKK}). However, we
will not use this pairing. 

A simple element $T_i-v$ acts now by  
\begin{equation}
\begin{aligned}
  M'=&\bigoplus_{\eps\in E} \Qv \eps,\\
(T_i-v) \eps&=\begin{cases} s_i \eps - v\eps, 
                  & (\eps_i, \eps_{i+1})=(+-),\\
                  0,
                  & (\eps_i, \eps_{i+1})=(--)\text{ or } (++),\\
               s_i\eps  -\vi \eps,
                  &  (\eps_i, \eps_{i+1})=(-+). 
          \end{cases}        
\end{aligned}
\end{equation}

Consider the space  $\P(k,n)$ of polynomials in
$x_1,\ldots,\, x_n$ of total degree $n-k$, and of degree  at most $1$ in
each $x_i$.  For any partition $\l$, denote by $x^{[\l]}$ the monomial
$w_\l(x_{k+1}\cdots x_n)$, the symmetric group acting now by permutation
of the $x_i$.  In other words, if $w_\l(\vac )= (\eps_1,\ldots,\,
\eps_n)$, then $x^{[\l]}$ is the product of the $x_i$'s for those $i$
such that $\eps_i=-$.

Consider the isomorphism of vector spaces 
\begin{equation}\lbb{isom}
\begin{aligned}
M'&\simeq \P(k,n)\\
w_\l(\vac )&\mapsto v^{-|\l| } x^{[\l]}.
\end{aligned}
\end{equation}
 Then $T_i -v$ induces the operator
$\nabla_i$, acting only on $x_i,x_{i+1}$ as follows:

\begin{equation}
  \begin{cases}
      \nabla_i(x_i)= v x_{i+1} -\vi x_i ,\\
      \nabla_i(1)=\nabla_i(x_ix_{i+1})= 0 ,\\
     \nabla_i(x_{i+1})= -v x_{i+1} +\vi x_i ,
          \end{cases}        
\end{equation}

Therefore $\nabla_i$ is the operator
$$
f\mapsto (v x_{i+1} -\vi x_i)\, \d_i(f)
$$
denoting by $\d_i$ the divided difference 
$$ 
f\mapsto  \frac{f- f^{s_i}}{x_i -x_{i+1}} 
$$
(for a more general action of the Hecke algebra on the ring of
polynomials, see \cite{LS2}, \cite{YB}).

We intend to show that divided differences easily furnish the \KL basis of
$\P(k,n)$ (i.e. the image of the \KL basis $C'_y, y\in W^J$ of $M'$). 

%\vskip 5mm
%{\tt If you fear that $P$ is confusing for the basis - not be confused
%with KL polynomials, take any letter which pleases you }
%
%\vskip 5mm
To any element $\eps:= w_\l(\vac )$ of $E$ one associates a polynomial 
$Q_\eps$ as follows

1) pair recursively  $-,+$ (as one pairs opening and closing parentheses)

2) replace each pair $(-,+)$, where $-$ is in position $i$ and $+$ in
position $j$, with a $ x_i -v^{j+1-i} x_j$

3) replace each single $-$, in position $i$, by $x_i$

  The product of all these factors by $v^{-|\l|}$, where
$|\l|= \l_1+\l_2+\cdots$,  is  by definition $Q_\eps$.

\begin{theorem}\lbb{newfactor}
Let $E$ be the set of sequences of $(+,-)$ of length $n$ with $k$ pluses.
Then the collection of polynomials $Q_\eps$, $\eps\in E$, is the
\KL basis of the space $\P(k,n)$. 
\end{theorem}

\begin{proof} 
 We shall show that
$$
Q_\eps = \nabla_j\cdots \nabla_h (x_1\cdots x_k)
$$
when $\eps = w_\l(\vac )$, and when $s_j\cdots s_h$ is a reduced
decomposition of $w_\l$.  Now, it is clear that the inverse image of
$Q_\eps$ in $M'$ is invariant under involution, and it is easy to
check the powers of $v$ to get that for $v=0$, it specializes to
$\eps$.

Assume by induction that we already know $Q_\eps$.
Let us add on the right of $\eps$ sufficiently many pluses, so that 
all minuses are now paired (the original polynomial is recovered from the
new one by specializing $x_{n+1}, x_{n+2}, \ldots$ to 0).
Take now any simple
transposition $s_i$ such that  $\eps_i=+, \eps_{i+1}= -$.
The variables $x_i, x_{i+1}$ involve two or one factor in $P_\eps$,
depending whether $\eps_i$ is paired or not. 
The only possible cases for those factors and their images under 
$\nabla_i$ are
\begin{align*}
 (x_{i-a} - v^{a+1} x_i)( x_{i+1} - v^{b+1} x_{i+b+1})
&\mapsto (x_{i-a} - v^{a+b+2} x_{i+b+1}) (\vi x_i-v x_{i+1}) \\
 ( x_{i+1} - v^{b+1} x_{i+b+1}) &\mapsto (\vi x_i-v x_{i+1})
\end{align*}
but now the 
new pairing of $-,+$ differs from the previous one
exactly in the places described by the  factors on the right.
\end{proof}

\begin{corollary}
Let $\si_j\cdots \si_h$ be a reduced decomposition of $w\in W^J$. Then
 the corresponding \KL element $C^{\prime J}_w\in M'$ is equal to 
$(T_j-v) \cdots  (T_h-v) (\vac)$. 
\end{corollary}

This factorization is equivalent to the one given in
\cite[Theorem 3.1]{FKK}. One can check on examples that this
factorization is compatible, via the duality between the two modules
$M$ and $M'$, with the factorization given by Theorem~\ref{main}.
However, deducing Theorem~\ref{main} from Theorem~\ref{newfactor} seems
more intricate than proving the two factorization properties directly.

%\medskip
%{\tt I have choosen an example with a non-trivial pairing, 
%which can be displayed in  matrix; example 1 is not more general, but
%too big}

\begin{example} Let $\l=[5,3,2]$ and $\mu=[5,3,3]$. Then  one has 

$$\begin{matrix}
{\s places} &1&2&3&4&5&6&7&8&9 \cr
\noalign{\smallskip \hbox to \hsize{\hrulefill}}
 {w_\l(\vac )}     & + &- &- &+           &-&+  &-&-         &+ \cr
\noalign{\smallskip \hbox to \hsize{\dotfill}}
                & + &- &  &            & &   &- &        &    \cr
{\s pairing}    &   &  &- &+           &-&+  &  &-       &+   \cr
\noalign{\smallskip \hbox to \hsize{\dotfill}}
 {\s polynomial} &   & x_2 &            & &   &  &x_7            \cr
            &   &    &(x_3 &-v^2x_4)&(x_5 &-v^2x_6)& &(x_8 &-v^2x_9)   \cr
\noalign{\smallskip \hbox to \hsize{\hrulefill}}
 {w_\mu(\vac )}  & + &- &- &-           &+&+  &-&-         &+ \cr
\noalign{\smallskip \hbox to \hsize{\dotfill}}
                & + &- &  &            & &   &- &        &    \cr
{\s pairing}    &   &  &- &            & &+  &  &        &    \cr
                &   &  &  &-           &+&   &  &-       &+  \cr
\noalign{\smallskip \hbox to \hsize{\dotfill}}
 {\s polynomial} &   & x_2 &            & &   &  &x_7            \cr
            &   &    &x_3          & & &-v^4x_6 &    & \cr
            &   &    & &(x_4       &-v^2x_5) && &(x_8 &-v^2x_9) \cr 
\noalign{\smallskip \hbox to \hsize{\hrulefill}}
\end{matrix} $$       

\medskip 
and thus 
\begin{equation}
\begin{aligned}
Q_{w_\l(\vac )} &=v^{-10}\, x_2 x_7\,(x_3-v^2x_4)\, 
                 (x_5-v^2x_6)\,(x_8-v^2x_9)         \cr
 Q_{w_\mu(\vac )} &=v^{-11}\, x_2 x_7\, (x_3-v^4x_6)\,
(x_4-v^2x_5)\,(x_8 -v^2x_9). 
\end{aligned}
\end{equation}
\end{example}

\medskip Note that the pairing between $-,+$, which was a key point in
the description of \KL polynomials for Grassmannians in \cite{LS}, is
provided by divided differences, starting from the monomial
$x_{k+1}\cdots x_n$.

%%%%%%%%%%%%%%%%%%%%%%%%%%%%%%%%%%%%%%%%%%%%%%%%%%%%%%%%%%%%%%%%%%%%%%
\bibliographystyle{amsalpha}

\end{document}